\documentclass[reqno]{amsart}
\usepackage{color}

\usepackage{cmbright}
\usepackage{amssymb,latexsym}
\usepackage[mathscr]{eucal}

\DeclareMathAlphabet{\cmrnrm}{OT1}{cmr}{m}{n}

\usepackage{graphicx}

\numberwithin{equation}{section}

\theoremstyle{plain}
\newtheorem{theorem}{Theorem}[section]
\newtheorem{lemma}[theorem]{Lemma}

\theoremstyle{remark}

\newtheorem{remarks}[theorem]{Remarks}
\newtheorem*{remarks*}{Remarks}
\newtheorem*{claim*}{Claim}
\newtheorem*{example*}{Example}

\newtheorem*{remark*}{Remark}
\newtheorem*{notation*}{Notation}

\newcommand{\figref}[1]{Figure \ref{#1}}
\newcommand{\lemref}[1]{Lemma \ref{#1}}

\newcommand{\secref}[1]{Section \ref{#1}}
\newcommand{\subsecref}[1]{Subsection \ref{#1}}
\newcommand{\subsubsecref}[1]{Subsubsection \ref{#1}}
\newcommand{\thmref}[1]{Theorem \ref{#1}}

\newenvironment{remlist*}{%
  \begin{remarks*}%
    \makebox{\phantom{a}}%
    \begin{enumerate}%
}{%
    \end{enumerate}%
  \end{remarks*}%
  }%
  
  \newenvironment{thmlist}{%
  \begin{enumerate}%
}{%
  \end{enumerate}%
  }%

\newcommand{\ls}[1]{
	\phantom{}_{\phantom{}_{#1}}
}

\newcommand{\Cut}{%
	\mathop{\mathrm{cut}}\nolimits%
}

\newcommand{\Det}{%
	\mathop{\mathrm{Det}}\nolimits%
}

\newcommand{\Diam}{%
	\mathop{\mathrm{diam}}\nolimits%
}

\newcommand{\Dist}{%
	\mathop{\mathrm{dist}}\nolimits%
}

\newcommand{\Sgn}{%
	\mathop{\mathrm{sgn}}\nolimits%
}

\newcommand{\Span}{%
	\mathop{\mathrm{span}}\nolimits%
}

\newcommand{\Vol}{%
	\mathop{\mathrm{Vol}}\nolimits%
}

\newcommand{\dVol}{%
	\mathop{\mathrm{dVol}}\nolimits%
}

\newcommand{\smeas}{{\mu^-}}
\newcommand{\tmeas}{{\mu^+}}

\newcommand{\sdom}{{\Omega^-}}
\newcommand{\tdom}{{\Omega^+}}

\newcommand{\bothdom}{{\Omega^\pm}}
\newcommand{\scost}{{c}}

\newcommand{\sfunc}{f}
\newcommand{\tfunc}{g}

\newcommand{\spot}{u}

\newcommand{\strans}{T}

\begin{document}


\title[Regularity of Potentials on Round Spheres]{On the Regularity of Optimal Transportation Potentials on Round Spheres}


\author[G.\,T. von~Nessi]{Greg~T. von~Nessi}
\address{Centre for Mathematics and Its Applications\\
	Australian National University\\
	Canberra ACT 0200, Australia}
\email{greg.vonnessi@anu.edu.au}
\urladdr{http://www.maths.anu.edu.au/~vonnessi/}



\begin{abstract}
In this paper the regularity of optimal transportation potentials defined on round spheres is investigated. Specifically, this research generalises the calculations done by Loeper, where he showed that the strong (A3) condition of Trudinger and Wang is satisfied on the round sphere, when the cost-function is the geodesic distance squared. In order to generalise Loeper's calculation to a broader class of cost-functions, the (A3) condition is reformulated via a stereographic projection that maps charts of the sphere into Euclidean space. This reformulation subsequently allows one to verify the (A3) condition for any case where the cost-fuction of the associated optimal transportation problem can be expressed as a function of the geodesic distance between points on a round sphere. With this, several examples of such cost-functions are then analysed to see whether or not they satisfy this (A3) condition.\end{abstract}

\maketitle

\section{Introduction}\label{in}
In the optimal transportation problem, one seeks to find an optimal mapping from one mass distribution to another, such that a given cost functional is minimised over all mappings which satisfy a change of variables criterion. Formally, one is given two topological spaces $\sdom$, $\tdom$, a cost-function $c:\sdom\times\tdom\to\mathbb{R}$ and two radon measures $\smeas$ and $\tmeas$ on $\sdom$ and $\tdom$ respectively, with the goal of finding an optimal map, $T:\sdom\to\tdom$, that minimises the functional
\begin{equation}
\int_\sdom c(x,T(x))\,d\sdom\label{in::eq1}
\end{equation}
over all mappings which satisfy the change of variables formula:
\begin{equation}
\int_{T^{-1}(\tdom)}(h\circ T)\,d\smeas=\int_\tdom h\,d\tmeas,\quad\forall h\in C^0(\tdom;\mathbb{R}).\label{in::eq2}
\end{equation}
$\tmeas$ is referred to as the \textit{push-forward} of $\smeas$ by the map $T$, which is notationally denoted as $T_\#\smeas=\tmeas$. The optimal transportation problem was first proposed by Monge in 1781 \cite{cit::mon1} and, in the above formulation, poses many obstacles to being analysed via standard methods in the direct calculus of variations (see \cite{cit::urb1,cit::von1} for details).\\

For the rest of the paper, it will be assumed that $d\smeas = f\,dx$ and $d\tmeas = g\,dy$. Also, for the time being, it will be assumed that both $\sdom$ and $\tdom$ are bounded subsets of $\mathbb{R}^n$. Even with these assumptions, none of the issues with Monge's original formulation become resolved. Indeed, it was not until Kantorovich \cite{cit::kan1,cit::kan2} introduced the dual formulation of Monge's original question in 1942, that any significant progress on this problem was made. Through Kantorovich's dual formulation it can be shown that an elliptic solution of the boundary value problem:
\addtocounter{equation}{1}
\begin{align}
\Det\left[\scost_{xx}(x,\strans_\spot(x))-D^2\spot(x)\right] &= \left|\Det\left[\scost_{xy}(x,\strans_\spot(x))\right]\right|\frac{\sfunc(x)}{\tfunc(\strans_\spot(x))},\quad x\in\sdom,\tag{\theequation a}\label{in::eq3a}\\
\strans_\spot(\sdom) &= \tdom\tag{\theequation b}\label{in::eq3b}.
\end{align}
yields a potential of Monge's original transportation problem, in a sense that if $T_u$ is subsequently defined as solving
\begin{equation}
Du=c_x(x,T_u(x))
\end{equation}
for a.e. $x\in\sdom$, then $T_u$ is indeed a minimiser of \eqref{in::eq1} satisfying the change of variables criterion embodied in \eqref{in::eq2} (See \cite{cit::eva2,cit::vil1,cit::von1} and references therein for details.) \eqref{in::eq3a} is an elliptic partial differential equation of a Monge-Amp\`ere type and is called the Optimal Transportation Equation. The boundary condition in \eqref{in::eq3b} is often referred to as a natural boundary condition.\\

Several conditions need to be placed on the cost-function in order to ensure existence and uniqueness of optimal maps. Specifically, denoting  $U\supseteq\sdom\times\tdom$, if the cost-function $c$ is real-valued in $C^4(U)$ and satisfies\\

\begin{enumerate}
\item[(A1)] For any $(x,y)\in U$ and $(p,q)\in D_x\scost(U)\times D_y\scost(U)$, there exists a unique $Y=Y(x,p)$, $X=X(q,y)$, such that $\scost_x(x,Y)=p$, $\scost_y(X,y)=q$.
\item[(A2)] For any $(x,y)\in U$,
\begin{equation*}
\Det\left[D^2_{xy}\scost\right]\neq 0,\label{pf::kc::eq1}
\end{equation*}
where $D_{xy}^2\scost$ is the matrix whose elements at the $i^\text{th}$ row and $j^{\text{th}}$ column is $\frac{\partial^2\scost}{\partial x_i\partial y_j}$,\\
\end{enumerate}

\noindent then existence and uniqueness (up to a constant) of optimal maps can be readily inferred from the proofs in \cite{cit::caf1,cit::gan_mcc1}. The regularity of optimal maps generated by solutions of the Optimal Transportation Equation for general cost-functions has only recently come into full understanding in the Euclidean space setting with the research presented in \cite{cit::ma_tru_wan1}, \cite{cit::tru_wan6}, \cite{cit::loe2} and \cite{cit::loe1}. Essentially, the cost-function, in addition to having (A1) and (A2) hold, needs to satisfy:\\

\begin{enumerate}
\item[(A3)] There exists a constant $C_0>0$ such that for any $(x,y)\in U$, and $\xi,\eta\in\mathbb{R}^n$ with $\xi\perp\eta$ such that
\begin{equation}
(\scost^{q,r}\scost_{ij,q}\scost_{r,st}-\scost_{ij,st})\scost^{s,k}\scost^{t,l}\xi_i\xi_j\eta_k\eta_l\le -C_0|\xi|^2|\eta|^2,\label{pf::kc::eq2}
\end{equation}
where $\scost_{i,j}(x,y)=\frac{\partial^2\scost(x,y)}{\partial x_i\partial y_j}$, and $[\scost^{i,j}]$ is the inverse matrix of $[\scost_{i,j}]$.\\
\end{enumerate}

\noindent in order to ensure that solutions of \eqref{in::eq3a} are locally smooth (provided $f$ and $g$ are smooth) \cite{cit::ma_tru_wan1}. The analogous global regularity result presented in \cite{cit::tru_wan6} only requires a degenerate form of the (A3) condition to hold that corresponds to when $C_0=0$. This degenerate version of the (A3) condition is labelled (A3w) in accordance with the notation in \cite{cit::tru_wan6}. While the research in \cite{cit::ma_tru_wan1,cit::tru_wan6} showed that the (A3) condition was sufficient for higher regularity of optimal transportation potentials, it was Loeper in \cite{cit::loe2} that showed that the (A3) condition was indeed necessary for these regularity results to hold. Specifically, Loeper showed that if (A3w) is violated, one can build a pair of $C^\infty$, strictly positive measures, supported on sets with the usual smoothness and convexity assumptions, so that the optimal potential is not even $C^1$; and thus, the corresponding optimal map is discontinuous. Before this work, it was not known whether the (A3w) condition was truly fundamental for potential function regularity or if it was simply a technical condition to make the \emph{a priori} estimates in \cite{cit::ma_tru_wan1} and \cite{cit::tru_wan6} work.

\begin{remark*}
It should be noted that in addition to the above criterion on the cost-fuction, additional convexity and smoothness conditions need to be satisfied by $\partial\bothdom$ in order for the global regularity estimates in \cite{cit::tru_wan6} to carry through; this will not be elaborated upon in this paper.
\end{remark*}

The author would like to thank Neil Trudinger and Xu-Jia Wang for their discussions and advisement regarding the research contained in this paper. 

\section[Some Properties of the (A3) Condition]{Some Properties of the (A3) Condition}\label{pf::kc::sp}
In order to carry out the stereographic reformulation later in the paper, some elementary properties of (A3) will be need. Specifically, the (A3) condition is symmetric in the $x$ and $y$ arguments, in addition to being invariant under coordinate transformations. For clarity, these properties are derived below.\\

To start off, it will be shown that the (A3) condition is symmetric in $x$ and $y$. To demonstrate this, $\tilde{\xi}$ is defined as
\begin{equation*}
\tilde{\xi}_k:=\scost_{q,k}\xi_q.
\end{equation*}
Using this definition to rewriting \eqref{pf::kc::eq2} with $\xi$ replaced by $\tilde{\xi}$, one sees that
\begin{equation}
(\scost^{p,q}\scost_{ij,p}\scost_{q,rs}-\scost_{ij,rs})\scost^{i,t}\scost^{j,h}\scost^{r,k}\scost^{s,l}\tilde{\xi}_t\tilde{\xi}_h\eta_k\eta_l\ge C_0|\xi|^2|\eta|^2,\label{pf::kc::sp::eq2}
\end{equation}
with a modified orthogonality criterion of
\begin{equation}
\eta_q\scost^{q,r}\tilde{\xi}_r=0.\label{pf::kc::sp::eq3}
\end{equation}
The symmetry of $x$ and $y$ in the (A3) criterion is now evident in the rewritten form embodied in \eqref{pf::kc::sp::eq2} and \eqref{pf::kc::sp::eq3}.
\begin{remark*}
The symmetry of \ref{pf::kc::sp::eq2} only makes sense when the mixed-Hessian of $\scost$ exists and will still hold if the constant $C_0$ degenerates. This will be important in the forth-coming calculations, where costs that depend on geodesic distance between points on a round sphere will be considered. Specifically, these costs have an ill-defined mixed-Hessian for argument pairs that are antipodal points of each other. This situation and its consequences on the results of this paper will be elaborated upon in Section \ref{sn::ac}.
\end{remark*}

Next, the expression depicted in \eqref{pf::kc::eq2} will be reduced to a less-cumbersome form, that also gives a better intuition as to what the (A3) condition itself actually means. Calculating, one see that
\begin{align}
D_{p_k}\scost_{ij}(x,Y(x,p)) &= \scost_{ij,q}\cdot D_{p_k}Y^q\notag\\
&= \scost_{ij,q}\scost^{q,k},\label{pf::kc::sp::eq4}
\end{align}
where the definition of $Y$ (as stated in (A1)) has been utilised to gain the second equality. Differentiation of \eqref{pf::kc::sp::eq4} subsequently yields
\begin{equation}
D^2_{p_lp_k}\scost_{ij}(x,Y(x,p))=\left(\scost_{ij,qr}\scost^{q,k}+\scost_{ij,q}\scost^{q,k}_{,r}\right)\scost^{r,l}.\label{pf::kc::sp::eq5}
\end{equation}
To proceed, a reduction of the term $\scost^{q,k}_{,r}$ is required. From previous notational definitions, it is understood that
\begin{equation*}
\scost_{i,q}c^{q,k}=\delta^k_i.
\end{equation*}
Differentiating this relation immediately produces the following relations:
\begin{align}
\scost_k^{i,j}(x,y) &= D_{x_k}\scost^{i,j}(x,y)\notag\\
&= -\scost^{i,q}\scost^{r,j}\scost_{kq,r}(x,y),\notag\\
\scost_{,k}^{i,j}(x,y) &= D_{y_k}\scost^{i,j}(x,y)\notag\\
&= -\scost^{i,q}\scost^{r,j}\scost_{q,kr}(x,y).\label{pf::kc::sp::eq7}
\end{align}
Combining \eqref{pf::kc::sp::eq5} with \eqref{pf::kc::sp::eq7}, it is now observed that
\begin{equation*}
D_{p_lp_k}\scost_{ij}(x,Y(x,p))=-(\scost^{q,r}\scost_{ij,q}\scost_{r,st}-\scost_{ij,st})\scost^{s,k}\scost^{t,l}.
\end{equation*}
Thus, the (A3) condition is equivalent to
\begin{equation}
D_{p_lp_k}\scost_{ij}(x,y)\xi_i\xi_j\eta_k\eta_l\le -C_0|\xi|^2|\eta|^2,\label{pf::kc::sp::eq9}
\end{equation}
for a positive constant $C_0$.\\

With \eqref{pf::kc::sp::eq9}, it is now a straight-forward calculation to verify that the (A3) condition is also invariant under coordinate transformations. Fixing $y$, consider an arbitrary change of coordinates in $x$ given by
\begin{equation*}
g(x)=x'.
\end{equation*}
An elementary calculation shows that
\begin{align}
\scost_{i}(x,y) &= [D_ig^q]\scost_{q}(x',y)\notag\\
\scost_{ij}(x,y) &= \left([D_ig^q][D_jg^r]\scost_{qr}(x',y)+[D_{ij}g^s]\scost_{s}(x',y)\right)\label{pf::kc::sp::eq11}
\end{align}
From (A1) and \eqref{pf::kc::sp::eq11}, it is observed that
\begin{align}
p_k' &= \scost_{k}(x',y)\notag\\
&= [D_kg^q]^{-1}\scost_{q}(x,y).\label{pf::kc::sp::eq12}
\end{align}
Using this in \eqref{pf::kc::sp::eq11} yields
\begin{equation}
\scost_{ij}(x',y)=\left([D_ig^q][D_jg^r]\scost_{qr}(x',y)+[D_{ij}g^s]p_s'\right).\label{pf::kc::sp::eq13}
\end{equation}
From here, one can use the chain-rule along with the relation in \eqref{pf::kc::sp::eq12} to deduce
\begin{equation}
D_{p_k}=[D_kg^q]^{-1}D_{p'_q}.\label{pf::kc::sp::eq14}
\end{equation}
Thus, from \eqref{pf::kc::sp::eq13} and \eqref{pf::kc::sp::eq14}, it has been shown that the left-hand side of \eqref{pf::kc::sp::eq9} is transformed to
\begin{equation*}
D_{p_lp_k}\scost_{ij}(x,y)[D_ig^q][D_jg^r][D_kg^s]^{-1}[D_lg^t]^{-1}\xi_q\xi_r\eta_s\eta_t.
\end{equation*}
Redefining $\xi_i$ as $[D_ig^q]\xi_q$ and $\eta_i$ as $[D_ig^q]^{-1}\eta_q$, finally shows that the (A3) condition is invariant under any change of coordinates, as the orthogonality criterion $\xi\perp\eta$ is preserved.\\

\section{Motivation}

In addition to showing the (A3w) condition to be sharp, in \cite{cit::loe2,cit::loe1} Loeper also studies the correlation between curvature and the Optimal Transportation Problem on Riemannian manifolds when the cost-function was of the form
\begin{equation}
\scost(x,y)=\frac{d^2(x,y)}{2},\label{sn::in::eq1}
\end{equation}
where $d(x,y)$ represents the geodesic distance between points $x,y\in\mathbb{M}^n$. One of the interesting aspects of this research was that Loeper verified the theory therein (in the case where $\mathbb{M}^n$ was taken to be a round sphere) by explicitly calculating that the cost-function depicted by \eqref{sn::in::eq1} satisfies the \emph{non}-degenerate form of the (A3) condition on round spheres. Specifically, it was shown that
\begin{equation}
D^2_{p_lp_k}\scost_{ij}(x,y)\xi_i\xi_j\eta_k\eta_l\to-\frac{2}{3},\quad\text{as}\ x\to y,\quad (x,y)\in\sdom\times\tdom\subset\mathbb{S}^n.\label{sn::in::eq2}
\end{equation}
This is particularly interesting as the cost-function in \eqref{sn::in::eq1} only satisfies the (A3w) condition in the Euclidean case but not the stronger (A3) condition (see \cite{cit::ma_tru_wan1,cit::tru_wan6}).\\

The disadvantage of Loeper's calculation in \cite{cit::loe1, cit::loe2} is that it relies heavily upon the specific geometric properties of the specific cost-function $\scost(x,y)=\frac{1}{2}d^2(x,y)$. In this paper, this calculation will be generalised on $\mathbb{S}^n$ to include cost-functions which are arbitrary functions of the geodesic distance with respect to the constant curvature metric on $\mathbb{S}^n$.

\section{Main Results}\label{sn::mr}
The main results of this paper centre around the verification of the (A3) criterion for various cost-functions (on round spheres) having the general form 
\begin{equation}
c(x,y)=f(d(x,y)),\label{sn::ac::eq1}
\end{equation}
where $d(x,y)$ represents the geodesic distance between $x$ and $y$ on the round sphere. The most important example thus considered is when $\scost(x,y)=\frac{1}{2}d^2(x,y)$. In this particular case, the following theorem (originally proven by Loeper and published in \cite{cit::loe1} and subsequently improved upon in \cite{cit::kim_mcc1,cit::fig_rif1,cit::fig_rif_vil1}) is verified:
\begin{theorem}[\cite{cit::loe1,cit::kim_mcc1,cit::fig_rif1,cit::fig_rif_vil1}]\label{sn::mr::thm1}
Given $\mathbb{S}^{n}$ be an embedded sphere in $\mathbb{R}^{n+1}$ with arbitrary radius, equipped with the round metric with an associated Riemannian geodesic distance $d$, if $\scost(x,y)=\frac{1}{2}d^2(x,y)$, then one has that $\scost$ satisfies the strong (A3) condition on $\mathbb{S}^{n+1}\times\mathbb{S}^{n+1}\setminus\{(x,-x)\,:\,x\in\mathbb{S}^n\}$.
\end{theorem}

The results in this paper extend the calculations in \cite{cit::loe1}. Indeed, the precise constant for which the (A3) condition is satisfied, for the situation where $\scost(x,y)=\frac{1}{2}d^2(x,y)$, will be calculated for round spheres of arbitrary radius. Loeper calculated this constant for the case where $R=1$; and the results in this paper indeed verify this.\\

Closely related to $\scost(x,y)=\frac{1}{2}d^2(x,y)$, is the cost-function of $$\scost(x,y)=2R^2\sin^2\left(\frac{d(x,y)}{2R}\right).$$ This scenario corresponds to the situation where $d$ is taken to be the chordal distance between two points in $\mathbb{S}^{n+1}$ relative to the ambient Euclidean space, into which the round sphere is embedded. For this case, the following analogous result to \thmref{sn::mr::thm1} will be proven:

\begin{theorem}\label{sn::mr::thm2}
Given $\mathbb{S}^{n}$ an embedded sphere in $\mathbb{R}^{n+1}$ with arbitrary radius, equipped with the round metric with an associated Riemannian geodesic distance $d$, if $\scost(x,y)=2R^2\sin^2\left(\frac{d(x,y)}{2R}\right)$, then one has that $\scost$ satisfies the strong (A3) condition on $\mathbb{S}^n\times\mathbb{S}^{n}\setminus\{(x,-x)\,|\,x\in\mathbb{S}^n\}$.
\end{theorem}

In addition to the two above situations, a few other examples of cost-functions are studied later in this paper.\\

It should be noted at this point that a cost-function satisfying the (A3) condition is not generally enough to guarantee higher regularity for Optimal Transportation potentials. Indeed, the Optimal Transportation Equation becomes singular for $y\in\Cut(x)$, for any cost-function depending on the geodesic distance between points $x$ and $y$ (specifically, $|\Det[D_{xy}^2\scost]|$ blows up). Thus, in order to have potential function regularity, bounds on transport vectors, that prevent optimal transport maps from mapping points to their cut-locus, are also needed in addition to the (A3) verification. This geometric criterion will have implications in the forthcoming analysis, which will be discussed in \secref{sn::ac}.

\section[Analysis of the (A3) Condition on $\mathbb{S}^n$]{Analysis of the (A3) Condition on $\mathbb{S}^n$}\label{sn::ac}
The following formulation relies on several key observations and simplifications that are unique in the specific case of analysing the (A3) condition on round spheres. First, the derivation of the Optimal Transportation Problem (and hence the Optimal Transportation Equation) only depends on the measure-space structure associated with $\sdom$ and $\tdom$. Indeed, the only way geometry can come into the formulation of the Optimal Transportation Problem is if the cost-function is defined as having some explicit geometric dependence. The cost-functions analysed in this paper all depend on the geodesic distance between points on $\mathbb{S}^n$, equipped with a round Riemannian metric. Thus, the underlying goal of the following calculation is to derive an explicit expression for the geodesic distance between two arbitrary, fixed points $x$ and $y$ lying on a round sphere. As will be discussed later, the round sphere is one of the very few manifolds (to the author's knowledge) where the geodesic distance between two arbitrary points can be explicitly represented. The rest of the simplifying observations will be stated as needed in the following formulation.\\

Ultimately, the analysis of the (A3) condition will be carried out in $\mathbb{R}^n$ (that is, in a local chart) with an associated, modified cost-function that yields an equivalent Optimal Transportation Problem, as compared to the one originally defined on $\mathbb{S}^n$. As mentioned earlier, it can not be assumed that Optimal Transportation maps do not move points to their cut-locus. Thus, in order to analyse the Optimal Transportation Problem on a manifold, via a single local chart, requires another geometric criterion or observation. If bounds on transport vectors exists such that optimal maps stay away from the cut-locus of particular points, then no additional criterion need be placed on the problem. As a broad class of cost-functions will be analysed in this paper, no such bounds have been proven for all the examples thus contained. Therefore, it will be assumed that $\Cut(\sdom)\cap\tdom=\emptyset$, in addition to the standard Optimal Transportation hypotheses placed on $\sdom$ and $\tdom$, unless stated otherwise.\\

Since only cost-functions will be considered that have the form depicted in \eqref{sn::ac::eq1}, it is possible to analyse the Optimal Transportation Problem in a local chart with a modified cost-function whose associated Optimal Transportation Problem in Euclidean space is equivalent to that of \eqref{sn::ac::eq1} on $\mathbb{S}^n$. Such a local chart and modified cost-function are formulated via a stereographic projection of $\mathbb{S}^n$ embedded in $\mathbb{R}^{n+1}$ onto the some arbitrary tangent space of $\mathbb{S}^n$. With this projection, an explicit expression for $d(x,y)$ can be derived in terms of $\hat{x}$ and $\hat{y}$: the projected coordinates on an arbitrary tangent space of $\mathbb{S}^n$.

\subsection{Stereographic Projection}\label{sn::ac::sp}%
Given that it is assumed that $\Cut(\sdom)\cap\tdom=\emptyset$, there will be no geometric issues with the forthcoming calculations being carried out in a single chart of $\mathbb{S}^n$. The possibility of relaxing this criterion will be discussed at the end of this paper, in \secref{sn::co}. As the following stereographic formulation is tantamount to analysing the (A3) condition through a specific coordinate transformation, it is recalled that the (A3) condition is invariant under coordinate transformation according to the calculations in \subsecref{pf::kc::sp}; and thus, the result of the following calculations will hold in general.

\subsubsection{The Half-Sphere Stereographic Projection}\label{sn::ac::sp::hs}
\begin{figure}[!hbt]
\centerline{\includegraphics[width=8.1cm]{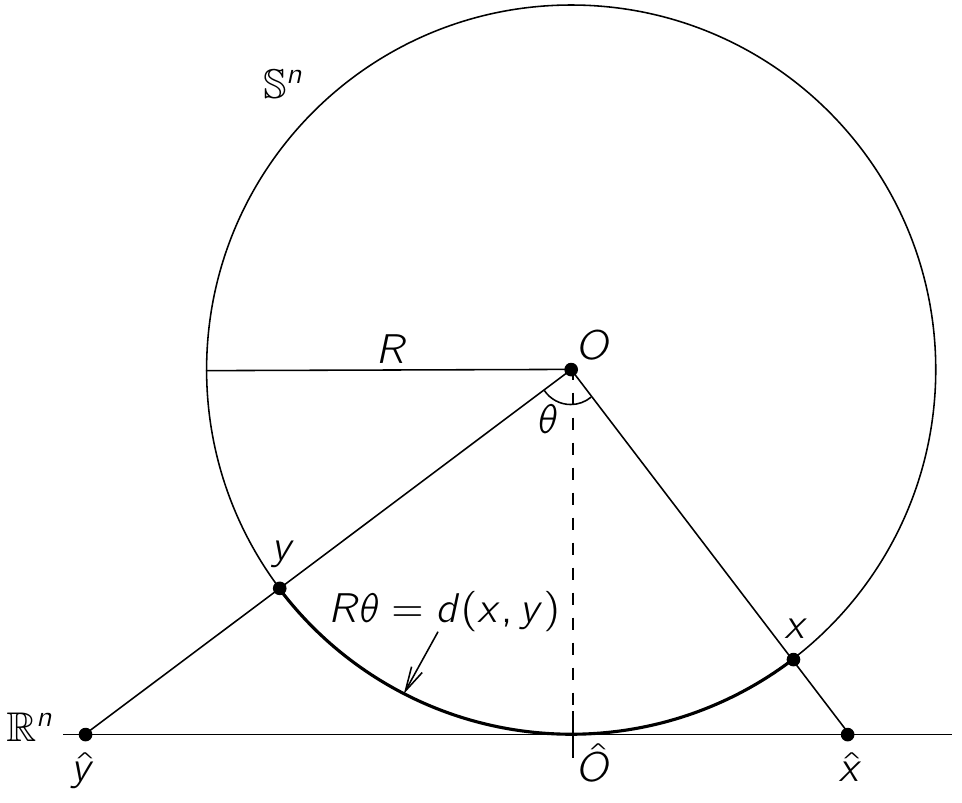}}
\caption{Stereographic projection for the half-sphere}\label{sn::ac::sp::hs::fig1}
\end{figure}

To begin, the modified form of \eqref{sn::ac::eq1} is first derived using the stereographic projection on the half-sphere depicted in \figref{sn::ac::sp::hs::fig1}. Following this derivation, it will then be described how half-sphere projection can be modified to project all points on the sphere to the plane, excluding the antipodal point of the intersection of the projective plane and sphere.\\

Utilising the ambient Euclidean geometry of $\mathbb{R}^{n+1}$, it is an elementary calculation that yields
\begin{align}
d(x,y) &= R\cdot\theta\notag\\
&= R\cdot\arccos\left(\frac{R^2+\hat{x}\cdot
\hat{y}}{\sqrt{R^2+|\hat{x}|^2}\sqrt{R^2+|\hat{y}|^2}}\right),\label{sn::ac::sp::hs::eq1}
\end{align}
where the origin of $\mathbb{R}^{n+1}$ is set to coincide with the centre of the sphere without any loss of generality.\\

From this calculation, the analysis of \eqref{sn::ac::eq1} on $\mathbb{S}^n$ is reduced to studying the Optimal Transportation Problem associated with the cost-function
\begin{equation}
\scost(\hat{x},\hat{y})=f\left(R\cdot\arccos\left(\frac{R^2+\hat{x}\cdot
\hat{y}}{\sqrt{R^2+|\hat{x}|^2}\sqrt{R^2+|\hat{y}|^2}}\right)\right)\label{sn::ac::sp::hs::eq2}
\end{equation}
between $\widehat{\sdom}$ and $\widehat{\tdom}$ on the local chart; and thus, \eqref{sn::ac::sp::hs::eq2} will now be analysed in the Euclidean setting.
\begin{remlist*}
\item Even though the ambient Euclidean geometry of $\mathbb{R}^{n+1}$ is used to derive \eqref{sn::ac::sp::hs::eq1}, having the sphere embedded in $\mathbb{R}^{n+1}$ is not technically required nor is the association of the local chart with a particular tanget plane of $\mathbb{S}^n$. Indeed, these notions are intuitive conveniences, in the sense that the equivalence of the Euclidean formulation of the Optimal Transportation Problem, to the original problem on $\mathbb{S}^n$, is immediate by this explicit stereographic projection.

\item A key point to the above calculation is that there is a simple relationship between $\theta$ and the geodesic distance between two \emph{arbitrary} points on $\mathbb{S}^n$. If one were to project rays from a point on $\mathbb{S}^n$ instead of its centre (as in the case of a full-sphere stereographic projection), there would be no simple relation between the angle of the projected rays and the associated geodesic distance between two arbitrary points. The strength of \eqref{sn::ac::sp::hs::eq1} is that it is independent of the particular choice of tangent plane on which the stereographic projection is performed.
\end{remlist*}

\subsubsection{The Full-Sphere Stereographic Projection}\label{sn::ac::sp::fs}
\begin{figure}[!hbt]
\centerline{\includegraphics[width=11cm]{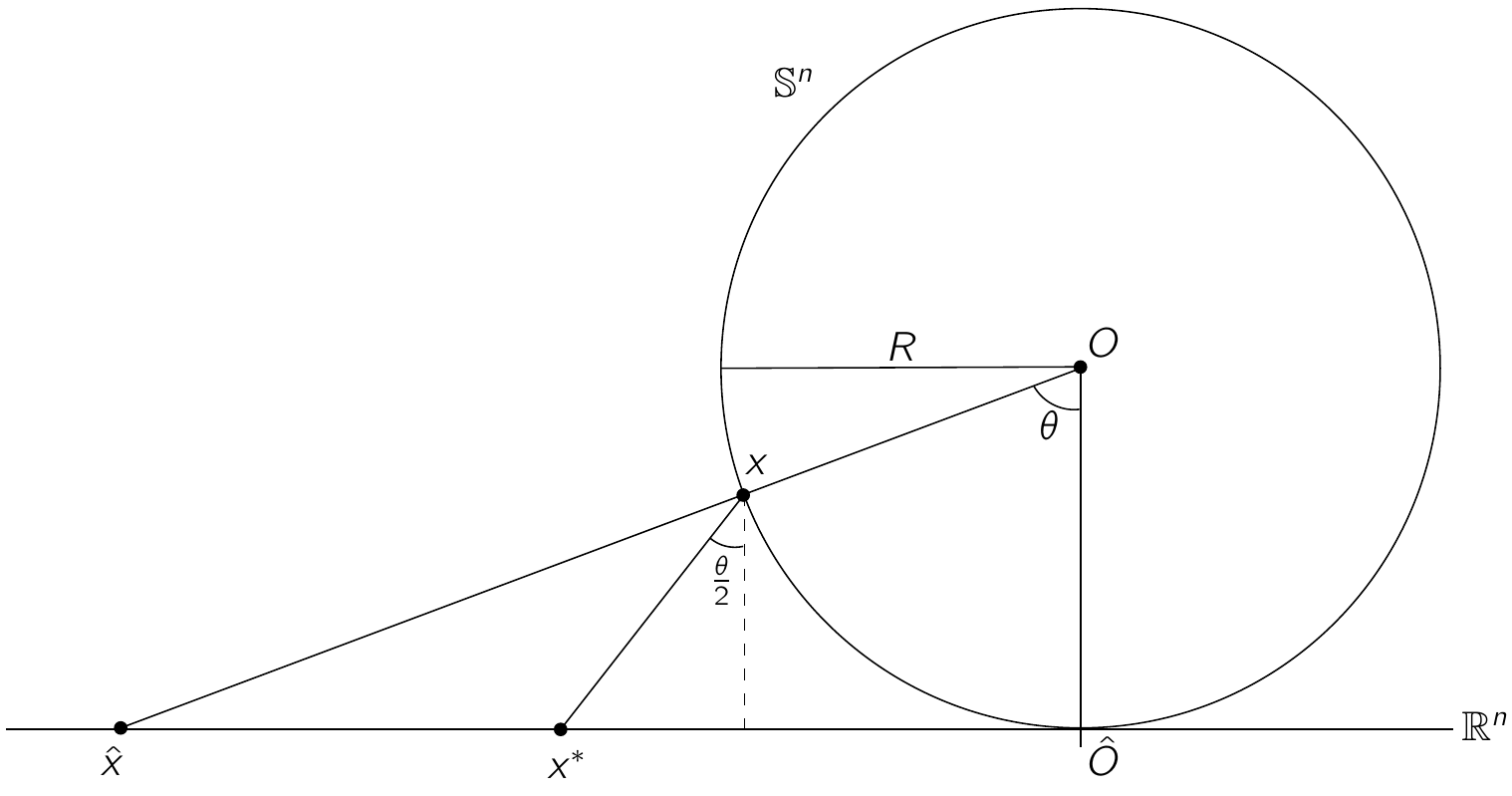}}
\caption{Stereographic projection for the full-sphere. Note that the projection ray corresponding to $\hat{x}$ has only been included for comparison to the half-sphere projection.}\label{sn::ac::sp::hs::fig1a}
\end{figure}

In order to analyse the (A3) condition between any points on a round sphere (excluding pairs of antipodal points), the half-sphere stereographic projection must be slightly modified. Specifically, the coordinate transform itself has to be changed such that the projected point is associated with a ray emitted from the centre of the sphere but subsequently has been refracted back toward the tangency point at the sphere's surface. This projection is depicted in \figref{sn::ac::sp::hs::fig1a}, where $x$ is taken to $x^*$ through this new coordinate transformation.\\

Using basic trigonometry, it is readily calculated that
\begin{equation}
|x^*|=2R\cdot\tan\left(\frac{\theta}{2}\right).\label{sn::ac::sp::fs::eq1}
\end{equation}
Moreover, one has that
\begin{equation}
\hat{x} = \left(1-\tan^2\left(\frac{\theta}{2}\right)\right)^{-1}x^*.\label{sn::ac::sp::fs::eq2}
\end{equation}
Combining the above two relations, yields
\begin{equation}
\hat{x}=\frac{4R^2}{4R^2-|x^*|^2}x^*.\label{sn::ac::sp::fs::eq3}
\end{equation}
This final relation allows one to rewrite \eqref{sn::ac::sp::hs::eq2} as
\begin{equation}
\scost(x^*,y^*)=f\left(R\cdot\arccos\left(\frac{1+\frac{16R^2}{(4R^2-|x^*|^2)(4R^2-|y^*|^2)}(x^*\cdot y^*)}{\sqrt{1+\frac{16R^2}{(4R^2-|x^*|^2)^2}|x^*|^2}\sqrt{1+\frac{16R^2}{(4R^2-|y^*|^2)^2}|y^*|^2}}\right)\right)\label{sn::ac::sp::fs::eq4}
\end{equation}\\

The advantage of \eqref{sn::ac::sp::fs::eq4} is that the relation holds between any two points on a sphere, excluding antipodal point pairs, for any arbitrarily chosen projective tangent plane. It is important to note that the construction \eqref{sn::ac::sp::fs::eq4} from the half-sphere stereographic projection is only done for expositional clarity; indeed, one could derive \eqref{sn::ac::sp::fs::eq4} with no use of the previous half-sphere projection.

\begin{notation*}
The forthcoming calculations will be based on the full-sphere relation in \eqref{sn::ac::sp::fs::eq4}, but $\hat{x}$ and $\hat{y}$ will be used instead $x^*$ and $y^*$ to indicate the full-sphere projected coordinates for the rest of the paper.
\end{notation*}

\subsubsection[Stereographic Reformulation of the (A3) Condition]{Stereographic Reformulation of the (A3) Condition}\label{sn::ac::sp::sr}%
With the relation \eqref{sn::ac::sp::fs::eq4}, it is now possible to formulate a new expression for the (A3) criterion in the current scenario. To do this, the geodesic distance on $\mathbb{S}^n$ is first differentiated with respect to the projected Euclidean coordinate. Subsequently, one may take the resulting expresion and set $\hat{x}$ to be the origin of the local chart without any loss of generality. This differentiation and origin selection ultimately give
\begin{equation}
d_{\hat{x}_i}(x,y)\Big|_{\hat{x}=0}=\frac{-\hat{y}_i}{|\hat{y}|}.\label{sn::ac::sp::sr::eq2}
\end{equation}
Differentiating twice with respect to the projected, Euclidean coordinates and subsequently choosing $\hat{x}$ to coincide with the origin yields
\begin{equation}
d_{\hat{x}_i\hat{x}_j}(x,y)\Big|_{\hat{x}=0} = -\frac{\hat{y}_i\hat{y}_j}{|\hat{y}|^3}+\frac{\delta_{ij}}{|\hat{y}|}.\label{sn::ac::sp::sr::eq4}
\end{equation}
Using these two relations, one can now readily differentiate \eqref{sn::ac::sp::fs::eq4} twice to find
\begin{equation}
\scost_{\hat{x}_i\hat{x_j}}=\frac{\hat{y}_i\hat{y}_j}{|\hat{y}|^2}\left(f''(d)-\frac{f'(d)}{|\hat{y}|}\right)+\delta_{ij}\frac{f'(d)}{|\hat{y}|},\label{sn::ac::sp::sr::eq5}
\end{equation}
where $\hat{x}=0$ has been assumed without loss of generality (since the choice was made \emph{after} differentiation).\\%

To proceed, both $\hat{y}$ and $d$ need to be represented in terms of the transportation vector, $\vec{p}$, which is defined in the formulation of the Optimal Transportation Problem by
\begin{equation*}
\vec{p}:=\nabla_{\hat{x}}\scost(\hat{x},\hat{y}),\label{sn::ac::sp::sr::eq6}
\end{equation*}
where $\scost(\hat{x},\hat{y})$ is recalled from \eqref{sn::ac::sp::fs::eq4}. Note that the vector notation on $p$ will be suppressed from this point forward without loss of clarity. From this expression and \eqref{sn::ac::sp::sr::eq2}, it is ascertained that
\begin{equation*}
p_i=-f'(d)\frac{\hat{y}_i}{|\hat{y}|};\label{sn::ac::sp::sr::eq7}
\end{equation*}
that is,
\begin{equation}
|p|=|f'(d)|,\quad\text{and}\quad\frac{p_i}{|p|}=-\frac{\hat{y}_i}{|\hat{y}|}.\label{sn::ac::sp::sr::eq8}
\end{equation}
On the other hand, it is clear from \figref{sn::ac::sp::hs::fig1} and \eqref{sn::ac::sp::fs::eq1}, with $\hat{x}=0$ being assumed, that one has
\begin{equation}
2R\cdot\tan\left(\frac{d}{2R}\right)=|\hat{y}|.\label{sn::ac::sp::sr::eq9}
\end{equation}
Moreover, since $\Cut(\sdom)\cap\tdom=\emptyset$ is assumed, it is guaranteed that $\frac{d}{R}<\pi$. 

\begin{remark*}
Again, it is important to note that the selection of the projective tangent plane is chosen and fixed only after the $x$ derivatives on the cost function are performed. Once the tangent plane is fixed, $\hat{x}$ can then be chosen as the origin in that plane.
\end{remark*}

Using both \eqref{sn::ac::sp::sr::eq8} and \eqref{sn::ac::sp::sr::eq9} \color{black}, one may rewrite \eqref{sn::ac::sp::sr::eq4} as
\begin{equation}
\scost_{\hat{x}_i\hat{x}_j}(p,d)=\frac{p_ip_j}{|p|^2}\left(f''(d)-f'(d)\cdot\mathcal{E}(d)\right)+\delta_{ij}f'(d)\cdot\mathcal{E}(d),\label{sn::ac::sp::sr::eq10}
\end{equation}
where
\begin{equation*}
\mathcal{E}(d;R):=\frac{1}{2R\cdot\tan\left(\frac{d}{2R}\right)}.
\end{equation*}
Defining the origin of the local chart to coincide with $\hat{x}$ effectively shifts all $p$ dependence onto the $\hat{y}$; and thus, differentiations with respect to $p$ variables may subsequently be applied directly to \eqref{sn::ac::sp::sr::eq10} with the understanding $\hat{x}=0$.\\

From \eqref{sn::ac::sp::sr::eq8} it is seen that $d$ does indeed have $p$ dependence. To differentiate $d$ with respect to $p$, one must first rewrite the expression in \eqref{sn::ac::sp::sr::eq8} as
\begin{equation*}
d=f'^{(-1)}\left(\Sgn(f')|p|\right).\label{sn::ac::sp::sr::eq11}
\end{equation*}
Differentiating implicitly, it is seen that
\begin{align}
d_{p_i} &= \frac{\Sgn(f')p_i}{f''|p|}\notag\\
&= \frac{p_i}{f'\cdot f''},\label{sn::ac::sp::sr::eq12}
\end{align}
where \eqref{sn::ac::sp::sr::eq8} has been used directly to ascertain the second equality, and the argument on $f$ has been suppressed without loss of clarity.\\

Now, the final form of the (A3) condition (associated with the general cost-function $f(d(x,y))$) is able to be written. Taking two arbitrary, unit vectors $\xi$ and $\eta$ such that $\xi\perp\eta$, one has that
\begin{align}
& D_{p_lp_k}^2\scost_{\hat{x}_i\hat{x}_j}(d)\xi_i\xi_j\eta_k\eta_l=\left(\frac{\mathcal{E}}{f'}+\frac{\mathcal{E}'}{f''}\right)+\frac{(p\cdot\xi)^2}{|p|^2}\left(\frac{f'''}{f'f''}-2\frac{f''}{f'^2}-\frac{\mathcal{E}'}{f''}+\frac{\mathcal{E}}{f'}\right)\notag\\
& \qquad\qquad+\frac{(p\cdot\eta)^2}{|p|^2}\left(\frac{\mathcal{E}'}{f''}-\frac{\mathcal{E}}{f'}+\frac{\mathcal{E}''f'}{f''^2}-\frac{\mathcal{E}'f'f'''}{f''^3}\right)+\frac{(p\cdot\xi)^2(p\cdot\eta)^2}{|p|^4}\left(\frac{f''''}{f''^2}\right.\notag\\
& \qquad\qquad\left.-5\frac{f'''}{f'f''}-\frac{f'''^2}{f''^3}+8\frac{f''}{f'^2}-\frac{\mathcal{E}''f'}{f''^2}+3\frac{\mathcal{E}'}{f''}+\frac{\mathcal{E}'f'f'''}{f''^3}-3\frac{\mathcal{E}}{f'}\right)\label{sn::ac::sp::sr::eq14}
\end{align}

\begin{remark*}
\eqref{sn::ac::sp::sr::eq14} gives the explicit, radial-scale dependency of the (A3) term in the current scenario. This has the interesting implication that it is possible to design cost-functions of the geodesic distance corresponding to the round sphere, so that for certain radii, the cost is strictly (A3) and for radii, it is not (A3) at all. This radial scale dependence is demonstrated explicitly in some of the forthcoming examples in \secref{sn::ex}.
\end{remark*}

For convenience, \eqref{sn::ac::sp::sr::eq14} is used to define four functions of $d$: $P_1(d)$, $P_2(d)$, $P_3(d)$ and $P_4(d)$, such that one may write
\begin{align}
D_{p_lp_k}^2\scost_{\hat{x}_i\hat{x}_j}(d)\xi_i\xi_j\eta_k\eta_l &= P_1(d)+\frac{(p\cdot\xi)^2}{|p|^2}P_2(d)+\frac{(p\cdot\eta)^2}{|p|^2}P_3(d)+\frac{(p\cdot\xi)^2(p\cdot\eta)^2}{|p|^4}P_4(d).\label{sn::ac::sp::sr::eq15}
\end{align}
To analyse \eqref{sn::ac::sp::sr::eq15}, various orientations of the vector $p$ relative to $\xi$ and $\eta$ must now be considered. These calculations can effectively be reduced down to analysing four cases.
\addtocounter{equation}{1}
\begin{description}
	\item[Case I:] If $p\perp\Span{(\xi,\eta)}$, then
		\begin{equation}
			D_{p_lp_k}^2\scost_{\hat{x}_i\hat{x}_j}(d)\xi_i\xi_j\eta_k\eta_l=P_1(d)=:O_1(d).\tag{\theequation a}\label{sn::ac::sp::sr::eq16a}
		\end{equation}
	\item[Case II:] If $p\,\|\,\xi$, then
		\begin{equation}
			D_{p_lp_k}^2\scost_{\hat{x}_i\hat{x}_j}(d)\xi_i\xi_j\eta_k\eta_l=P_1(d)+P_2(d)=:O_2(d).\tag{\theequation b}\label{sn::ac::sp::sr::eq16b}
		\end{equation}
	\item[Case III:] If $p\,\|\,\eta$, then
		\begin{equation}
			D_{p_lp_k}^2\scost_{\hat{x}_i\hat{x}_j}(d)\xi_i\xi_j\eta_k\eta_l=P_1(d)+P_3(d)=:O_3(d).\tag{\theequation c}\label{sn::ac::sp::sr::eq16c}
		\end{equation}
	\item[Case IV:] If $p\in\Span{(\xi,\eta)}$ and $p\cdot\xi=p\cdot\eta$, then
		\begin{equation}
			D_{p_lp_k}^2\scost_{\hat{x}_i\hat{x}_j}(d)\xi_i\xi_j\eta_k\eta_l=P_1(d)+\frac{P_2(d)}{2}+\frac{P_3(d)}{2}+\frac{P_4(d)}{4}=:O_4(d).\tag{\theequation d}\label{sn::ac::sp::sr::eq16d}
		\end{equation}
\end{description}
Here, four new functions of $d$ are again defined for clarity. Subsequently, the combined negativity of $O_1$, $O_2$, $O_3$ and $O_4$ is tantamount to the (A3) condition being satisfied. This is formally stated in the following lemma:
\begin{lemma}\label{sn::ac::sp::sr::lem1}
Given a cost-function of the form $\scost(x,y)=f(d(x,y))$ and $\mathbb{S}^{n}$ an embedded sphere in $\mathbb{R}^{n+1}$ with arbitrary radius $R$, equipped with the round Riemannian metric with an associated geodesic distance $d$, then the following statements are true:
\begin{thmlist}
\item If $O_i(d)\le C<0$ for $i\in\{2,3,4\}$ for all {\color{black}$d\in\left(0,R\pi\right)$} and $n>2$, then the strong (A3) condition is satisfied with constant $C$, for the cost-function.
\item If $O_i(d)<0$ for $i\in\{2,3,4\}$ for all {\color{black}$d\in\left(0,R\pi\right)$} and $n>2$, then the (A3w) condition is satisfied for the cost-function.
\item If $O_i(d)\le C<0$ for $i\in\{1,2,3,4\}$ for all  {\color{black}$d\in\left(0,R\pi\right)$} and $n\ge 2$, then the strong (A3) condition is satisfied with constant $C$ for the cost-function.
\item If $O_i(d)<0$ for $i\in\{1,2,3,4\}$ for all  {\color{black}$d\in\left(0,R\pi\right)$} and $n\ge 2$, then the (A3w) condition is satisfied for the cost-function for $n\ge2$.
\end{thmlist}
\end{lemma}

{\color{black}With the above calculations and recalling the statements made after \eqref{sn::ac::sp::sr::eq9}}, verification of the (A3w) condition is now reduced to verifying that \eqref{sn::ac::sp::sr::eq16a}--\eqref{sn::ac::sp::sr::eq16d} are non-positive for $d\in\left(0,R\pi\right)$. Correspondingly, \eqref{sn::ac::sp::sr::eq16a}--\eqref{sn::ac::sp::sr::eq16d} being strictly negative indicate that the strong (A3) condition is satisfied. The primary benefit gained by the use of \eqref{sn::ac::sp::sr::eq16a}--\eqref{sn::ac::sp::sr::eq16d}, is that one may uses these relations to easily generate explicit, analytic expressions that allow for straight-forward verification of the (A3) criterion for a general class of $f$. Such calculations will be the focus of the next section.

\section{Examples}\label{sn::ex}
In this section, various examples of cost-functions will be analysed using \eqref{sn::ac::sp::sr::eq16a}--\eqref{sn::ac::sp::sr::eq16d} to see if they do indeed satisfy the (A3w) or the (A3) condition. These examples include the some of the cost-functions analysed in \cite{cit::tru_wan6} and \cite{cit::ma_tru_wan1} for the Euclidean case and are encompassed by the general cost-function
\begin{equation}
\scost(x,y)=f(d(x,y)),\label{sn::ex::eq1}
\end{equation} 
where $d(x,y)$ represents the geodesic distance between $x$ and $y$ with respect to the underlying Riemannian manifold.
\begin{remark*}
The forthcoming results will be conveyed with the understanding that they hold for all $d\in\left(0,R\pi\right)$, unless otherwise indicated.
\end{remark*}

\subsection{$\scost(x,y)=\frac{1}{2}d^2(x,y)$}\label{sn::ex::d2}%
The first example consider is the cost-function that has been the most studied out of all the possibilities encompassed by \eqref{sn::ex::eq1}: $\scost(x,y)=\frac{1}{2}d^2(x,y)$. Indeed, this is the only cost-function for which the stereographic formulation may be applied without any extraneous geometric conditions, that serve to validate analysing the Optimal Transportation Problem in only one chart. This is due to the gradient estimate first proved in \cite{cit::mcc1} for compact Riemannian manifolds, which was then improved upon in the case of round sphere in \cite{cit::del_loe1}. To the author's knowledge, there exists no such gradient bounds for the general cost-function depicted in \eqref{sn::ex::eq1}. This estimate will now be briefly reviewed, as it will also motivate points of the discussion at the end of this paper.\\

As was stated above, the scenario where $\scost(x,y)=\frac{1}{2}d^2(x,y)$ has been studied previously in \cite{cit::mcc1} in the context of Optimal Transportation on general Riemannian manifolds; and it has been specifically analysed on the round sphere in both \cite{cit::loe1} and \cite{cit::del_loe1}. In \cite{cit::mcc1}, McCann studies gradient mappings defined on Riemannian manifolds, which are mappings of the form
\begin{equation*}
G_\phi(m):=\exp_m(\nabla_m\phi),\label{sn::ex::d2::eq1}
\end{equation*}
where $\phi$ is the associated gradient-potential of the mapping. McCann shows that such maps are indeed optimal in the transportation of measures on Riemannian manifolds for the cost-function $\scost(x,y)=\frac{1}{2}d^2(x,y)$. In addition to this, McCann also shows that if $\phi$ is $\scost$-convex, then the length of it's gradient can not exceed the diameter of the manifold. Translating into the context of stereographic projections of round spheres, this means that the transport vector defined in \eqref{sn::ac::sp::sr::eq5} is not necessarily bounded in the projected local chart, but it is well-defined there. Delano\"e and Loeper improved this result on $\mathbb{S}^n$ by proving the following gradient bound in \cite{cit::del_loe1}:
\begin{theorem}[\cite{cit::del_loe1}]\label{sn::ex::d2::thm1}
Given $\phi:\mathbb{S}^n\mapsto\mathbb{R}$ a $\scost$-convex function, such that
\begin{equation}
\int_{G_\phi^{-1}(\mathbb{S}^n)}(h\circ G_\phi)\,\dVol=\int_{\mathbb{S}^n}(h\cdot\rho)\,\dVol\label{sn::ex::d2::eq2}
\end{equation}
for some $\rho\in L^\infty(\mathbb{S}^n,\dVol)$ with any $h\in C^0(\mathbb{S}^n)$, where $\dVol$ stands for the canonical Lebesgue measure on $\mathbb{S}^n$, then the following estimate holds a.e.
\begin{equation}
|d\phi|\le\pi-\frac{1}{2\pi}\left\{\frac{1}{\|\rho\|\ls{L^\infty(\mathbb{S}^n)}}\left[\frac{n\Vol(\mathbb{S}^n)}{2\Vol(\mathbb{S}^{n-1})}\right]^2\right\}^{1/n}\label{sn::ex::d2::eq3}.
\end{equation}
\end{theorem}
\begin{remark*}
Recalling that the notation remarked upon in \secref{in}, the expression in \eqref{sn::ex::d2::eq2} is often denoted as
\begin{equation*}
{G_\phi}_\#\dVol=\rho\dVol.\label{sn::ex::d2::eq4}
\end{equation*}
\end{remark*}

\thmref{sn::ex::d2::thm1} correlates to the transportation vector $p$, as defined in \eqref{sn::ac::sp::sr::eq5}, being strictly bounded in the stereographically projected local chart. Thus, the assumption that $\Cut(\sdom)\cap\tdom=\emptyset$ is not needed in the case where $\scost(x,y)=\frac{1}{2}d^2(x,y)$. The possibility of extending \eqref{sn::ex::d2::eq3} to the general case of $\scost(x,y)=f(d(x,y))$, will be discussed at the end of this paper. The (A3) condition will now be analysed for the case when $\scost(x,y)=\frac{1}{2}d^2(x,y)$, using the results of \secref{sn::ac}.\\

\begin{figure}[!htb]
\centerline{\includegraphics[width=\textwidth]{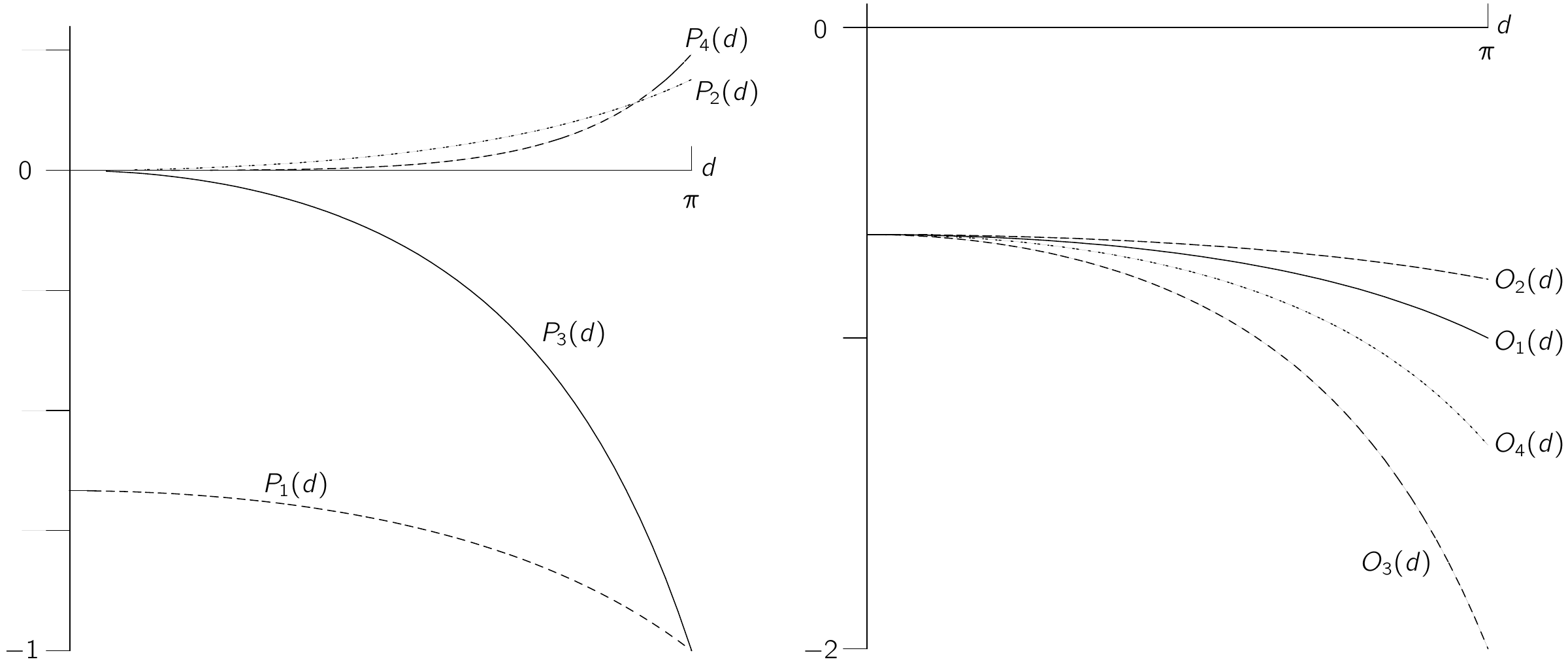}}
\caption{$P_i(d)$ and $O_i(d)$ corresponding to $\scost(x,y)=\frac{1}{2}d^2(x,y)$ with $R=1$}\label{sn::ex::d2::fig1}
\end{figure}

Using \eqref{sn::ac::sp::sr::eq14}, it is calculated that
\addtocounter{equation}{1}
\begin{align}
P_1(d) &= \frac{\cos\left((2R)^{-1}d\right)}{2Rd\sin\left((2R)^{-1}d\right)}-\frac{1}{4R^2\sin^2\left((2R)^{-1}d\right)}\tag{\theequation a}\label{sn::ex::d2::eq5a}\\
P_2(d) &= -\frac{2}{d^2}+\frac{\cos\left((2R)^{-1}d\right)}{2Rd\sin\left((2R)^{-1}d\right)}+\frac{1}{4R^2\sin^2\left((2R)^{-1}d\right)}\tag{\theequation b}\label{sn::ex::d2::eq5b}\\
P_3(d) &= -\frac{\cos\left((2R)^{-1}d\right)}{2Rd\sin\left((2R)^{-1}d\right)}-\frac{1}{4R^2\sin^2\left((2R)^{-1}d\right)}+2\frac{d\cos\left((2R)^{-1}d\right)}{8R^3\sin^3\left((2R)^{-1}d\right)}\tag{\theequation c}\label{sn::ex::d2::eq5c}\\
P_4(d) &= \frac{8}{d^2}-2\frac{d\cos\left((2R)^{-1}d\right)}{8R^3\sin^3\left((2R)^{-1}d\right)}-3\frac{\cos\left((2R)^{-1}d\right)}{2Rd\sin\left((2R)^{-1}d\right)}-3\frac{1}{4R^2\sin^2\left((2R)^{-1}d\right)}.\tag{\theequation d}\label{sn::ex::d2::eq5d}
\end{align}
As stated in \subsecref{sn::ac::sp}, \eqref{sn::ex::d2::eq5a}--\eqref{sn::ex::d2::eq5d} are considered for $d\in\left(0,R\pi\right)$. First, the limits as $d\to 0$ are calculated to be
\begin{equation*}
\lim_{d\to0}P_1(d)=-\frac{2}{3R^2},\ \lim_{d\to0}P_2(d)=0,\ \lim_{d\to0}P_3(d)=0,\ \lim_{d\to0}P_4(d)=0.\label{sn::ex::d2::eq6}
\end{equation*}
Using these limits to analyse $O_1(d)$, $O_2(d)$, $O_3(d)$ and $O_4(d)$, it follows that
\begin{equation*}
\lim_{d\to 0}D_{p_lp_k}^2\scost_{\hat{x}_i\hat{x}_j}(d)\xi_i\xi_j\eta_k\eta_l=-\frac{2}{3R^2},\quad\forall p\in\mathbb{R}^n.\label{sn::ex::d2::eq7}
\end{equation*}
This confirms the result of Loeper presented at the end of \cite{cit::loe1} for the case where $R=1$. Using \eqref{sn::ex::d2::eq5a}--\eqref{sn::ex::d2::eq5d}, it is also calculated that
\begin{equation*}
O'_i(d)<0,\quad\forall d\in\left(0,R\pi\right),\ \forall R>0,\label{sn::ex::d2::eq8}
\end{equation*}
where $i\in\{1,2,3,4\}$. \figref{sn::ex::d2::fig1} below contains plots of both sets of functions $P_i(d)$ and $O_i(d)$, in the case when $R=1$.\\

With the above calculations, it has been shown that the (A3) condition is indeed satisfied on the round sphere of any radius, for the cost of $\scost(x,y)=\frac{1}{2}d^2(x,y)$. As noted earlier, this is in contrast to the Euclidean case where the cost-function $\scost(x,y)=\frac{1}{2}|x-y|^2$ only satisfies the (A3w) condition, as the (A3) term goes to $0$ as $y\to x$ (see \cite{cit::ma_tru_wan1,cit::tru_wan6}). This contrast between the Euclidean and spherical cases indicates that the underlying geometry (manifested through the cost-function) does indeed affect the local regularity of the Optimal Transportation Problem.\\

\subsection{$\scost(x,y)=2R^2\sin^2\left(\frac{d(x,y)}{2R}\right)$}\label{sn::ex::s2}%
\begin{figure}[!hbt]
\centerline{\includegraphics[width=6.38cm]{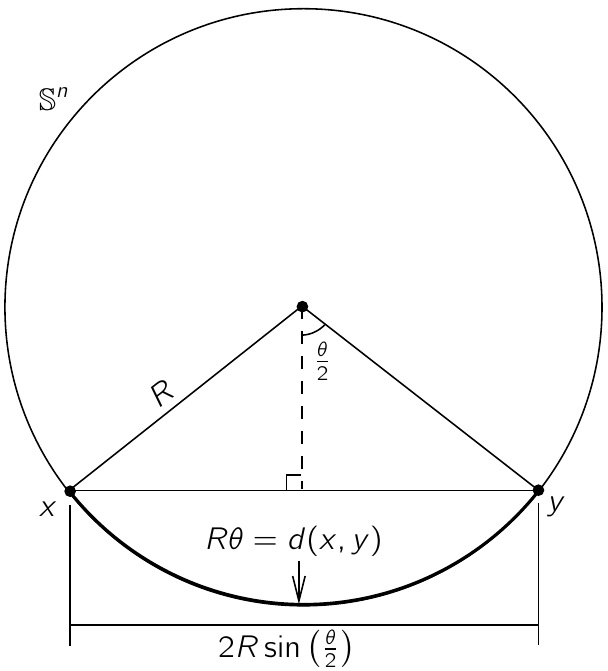}}
\caption{Chordal distance between points}\label{sn::ex::s2::fig1}
\end{figure}

If one were to consider the round sphere embedded in $\mathbb{R}^{n+1}$ equipped with a Euclidean metric, the cost-function of $2R^2\sin^2\left(\frac{1}{2R}d(x,y)\right)$ is equivalent to the example studied in \subsecref{sn::ex::d2} with $d(x,y)$ taken to being the geodesic distance of the Euclidean space into which the sphere is embedded. This situation is depicted in \figref{sn::ex::s2::fig1} on the next page.\\

Again, using the results of \subsubsecref{sn::ac::sp::sr}, it is calculated that
\addtocounter{equation}{1}
\begin{align}
P_1(d) &= \frac{1}{4R^2\left[1-2\cos^2\left((2R)^{-1}d\right)\right]}\tag{\theequation a}\label{sn::ex::s2::eq1a}\\
P_2(d) &= \frac{4\sin^2\left((2R)^{-1}d\right)\cos^2\left((2R)^{-1}d\right)}{4R^2\left[1-8\cos^6\left((2R)^{-1}d\right)+12\cos^4\left((2R)^{-1}d\right)-6\cos^2\left((2R)^{-1}d\right)\right]}\tag{\theequation b}\label{sn::ex::s2::eq1b}\\
P_3(d) &= 0\tag{\theequation c}\label{sn::ex::s2::eq1c}\\
P_4(d) &= \frac{2\sin^4\left((2R)^{-1}d\right)\cos^2\left((2R)^{-1}d\right)}{4R^2\left[8\cos^6\left((2R)^{-1}d\right)-12\cos^4\left((2R)^{-1}d\right)+6\cos^2\left((2R)^{-1}d\right)-1\right]}.\tag{\theequation d}\label{sn::ex::s2::eq1d}
\end{align}
From these relations, the following limits are calculated:
\begin{equation*}
\lim_{d\to0}P_1(d)=-\frac{1}{R^2},\ \lim_{d\to0}P_2(d)=0,\ \lim_{d\to0}P_3(d)=0,\ \lim_{d\to0}P_4(d)=0;
\end{equation*}
thus,
\begin{equation}
\lim_{d\to 0}D^2_{p_lp_k}\scost_{\hat{x}_i\hat{x}_j}(d)\xi_i\xi_j\eta_k\eta_l=-\frac{1}{R^2},\quad\forall p\in\mathbb{R}^n.\label{sn::ex::s2::eq3}
\end{equation}
Finally, an elementary calculation using \eqref{sn::ex::s2::eq1a}--\eqref{sn::ex::s2::eq1d} in the definitions of $O_i(d)$ yields 
\begin{equation}
O'_i(d)<0\quad\forall d\in\left(0,R\pi\right)\quad\forall R>0.\label{sn::ex::s2::eq4}
\end{equation}

Given the properties depicted in \eqref{sn::ex::s2::eq3} and \eqref{sn::ex::s2::eq4}, it follows that this cost-function does indeed satisfy the strong (A3) condition on the half-sphere. This is readily confirmed by the plots of $P_i(d)$ and $O_i(d)$ in the case where $R=1$, shown in \figref{sn::ex::s2::fig2} below.

\begin{figure}[!hbt]
\centerline{\includegraphics[width=\textwidth]{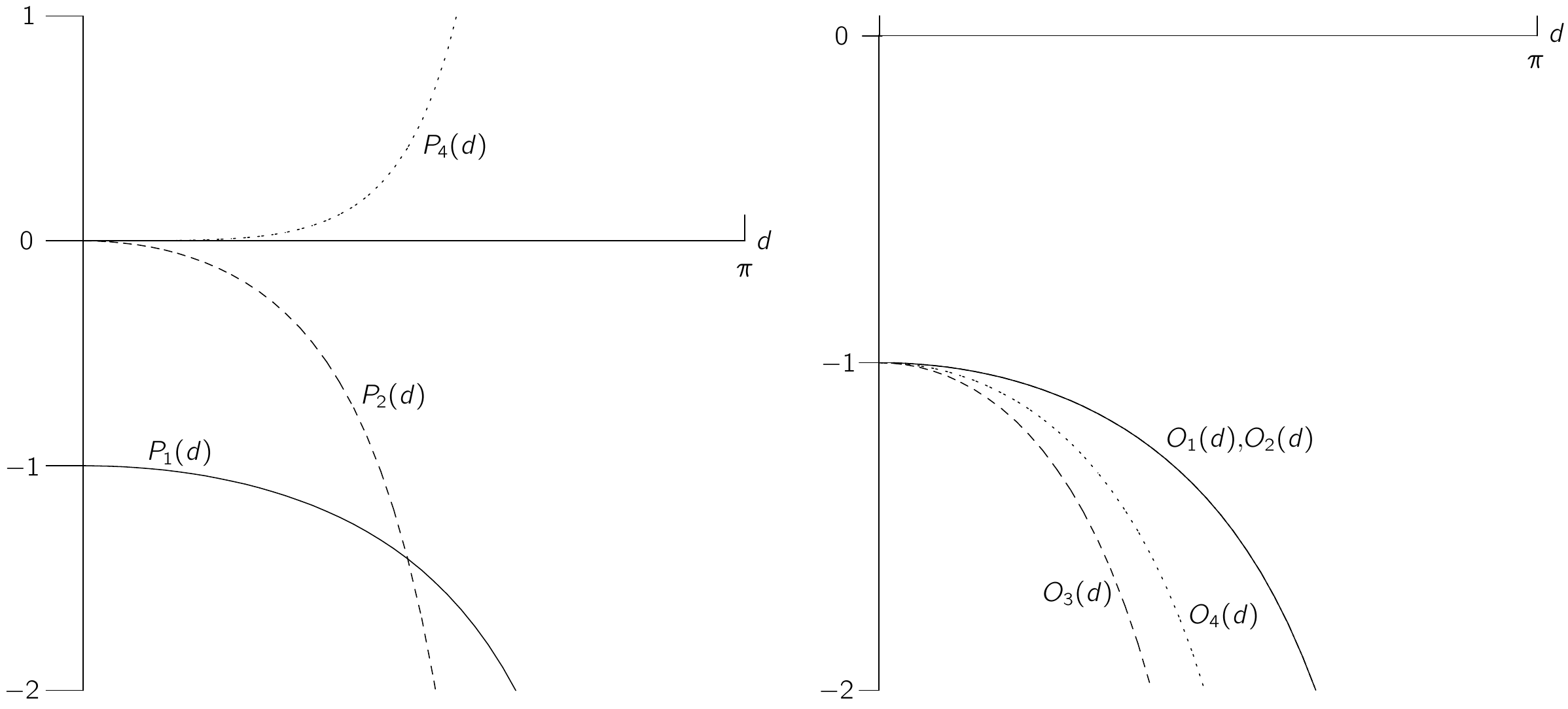}}
\caption{$P_i(d)$ and $O_i(d)$ corresponding to $\scost(x,y)=2R^2\sin^2\left(\frac{d(x,y)}{2R}\right)$ with $R=1$}\label{sn::ex::s2::fig2}
\end{figure}

\subsection{$\scost(x,y)=\sqrt{1-d^2(x,y)}$}\label{sn::ex::pm}
This example gives a demonstration as to just how important $R$ scaling can be in affecting the (A3) condition. This radial scaling dependence manifests itself through the limiting behaviour of the orientation terms as $d\to 0$. Using the general expression for $O_i(d)$, one readily calculates that
\begin{equation*}
\lim_{d\to 0}O_i(d)=\frac{2}{3R^2}-1\quad\text{for}\ i\in\{1,2,3,4\}.
\end{equation*}

Further calculation shows that for $R>\sqrt{2/3}$, and $\Diam(\sdom\cup\tdom)<1$, the (A3) condition is satisfied. For $R=\sqrt{2/3}$, the (A3) condition is satisfied if $\Dist(\sdom,\tdom)>0$ and $\Diam(\sdom\cup\tdom)<1$. If one only has that $\Dist(\sdom,\tdom)\ge0$ and $\Diam(\sdom\cup\tdom)\le1$, then only the (A3w) condition is satisfied. In other cases where $R\ge\sqrt{2/3}$, the (A3) condition will be violated.\\

If $R<\sqrt{2/3}$ and $\Diam(\sdom\cup\tdom)<1$ with $\Dist(\sdom,\tdom)>h^*$, where $h^*>0$ solves $P_1(h^*)=0$, then $\scost(x,y)=\sqrt{1-d^2}$ satisfies the strong (A3) condition. This cost-function will violate the (A3) condition in all other cases where $R<\sqrt{2/3}$.

\subsection{$\scost(x,y)=-\sqrt{1-d^2(x,y)}$}\label{sn::ex::mm}
This example is closely related to the previous example, in that it also has a dependence on the radial scaling of the round sphere. In particular, one has that
\begin{equation*}
\lim_{d\to 0}O_i(d)=1-\frac{2}{3R^2}\quad\text{for}\ i\in\{1,2,3,4\}.
\end{equation*}
From this, it is observed that for $R\ge\sqrt{2/3}$, the (A3) condition will be violated at least for some small values of $d$. Indeed, further calculation indicates that
\begin{equation*}
O_i(d)\ge0\quad\text{for}\ i\in\{1,2,3\},\ \forall d\in[0,1],
\end{equation*}
for any $R\ge\sqrt{2/3}$, thus proving the (A3) condition will be violated in this case.\\

If $R<\sqrt{2/3}$ and $\Diam(\sdom\cup\tdom)<\min\{h^*,R,1\}$, where $h^*>0$ solves $P_2(h^*)=0$, then the strong (A3) condition will be satisfied. In all other cases where $R<\sqrt{2/3}$, the (A3) condition will be violated.

\subsection{$\scost(x,y)=\sqrt{1+d^2(x,y)}$}\label{sn::ex::pp}
This example bears some resemblance to the previous one; but in this scenario, there is no sign-symmetry broken by a variance in the radial scaling. Indeed, it is calculated that
\begin{equation*}
\lim_{d\to 0}O_i(d)=-\frac{2}{3R^2}-1,\quad\text{for}\ i\in\{1,2,3,4\},
\end{equation*}
along with 
\begin{equation*}
O_i'(d)\le0,\quad i\in\{1,2,3,4\}.
\end{equation*}
In this scenario, the (A3) condition will be strongly satisfied for all values of $R>0$.

\subsection{$\scost(x,y)=-\sqrt{1+d^2(x,y)}$}\label{sn::ex::mp}
In this example, it is calculated that
\begin{equation*}
O_i(d)>0,\quad i\in\{1,2,3,4\},
\end{equation*}
for any $R>0$. Thus, the (A3) condition will always be violated for this cost-function.

\subsection{$\scost(x,y)=\pm\frac{1}{m}d^m(x,y)$}\label{sn::ex::dm}
In this subsection, a generalisation of the example in \subsecref{sn::ex::d2} is considered:
\begin{equation}
\scost(x,y)=\pm\frac{d^m}{m}(x,y)\ m\neq 0,\quad \scost(x,y)=\pm\log(d(x,y))\ m=0.\label{sn::ex::dm::eq1}
\end{equation}
Calculating as before, one finds the following:
\addtocounter{equation}{1}
\begin{align}
O_1(d) &= \pm\frac{(m-1)2R\sin\left((2R)^{-1}d\right)\cos\left((2R)^{-1}d\right)-d}{(m-1)4R^2d^{m-1}\sin^2\left((2R)^{-1}d\right)}\tag{\theequation a}\label{sn::ex::dm::eq2a}\\
O_2(d) &= \mp\frac{m2R\sin\left((2R)^{-1}d\right)-2d\cos\left((2R)^{-1}d\right)}{2Rd^m\sin\left((2R)^{-1}d\right)}\tag{\theequation b}\label{sn::ex::dm::eq2b}\\
O_3(d) &= \mp\frac{m2R\sin\left((2R)^{-1}d\right)-2d\cos\left((2R)^{-1}d\right)}{(m-1)^28R^3d^{m-2}\sin^3\left((2R)^{-1}d\right)}\tag{\theequation c}\label{sn::ex::dm::eq2c}\\
O_4(d) &= \pm\frac{\cos\left((2R)^{-1}d\right)}{2(m-1)^28R^3d^{m-3}\sin^3\left((2R)^{-1}d\right)}\pm\frac{\cos\left((2R)^{-1}d\right)}{8Rd^{m-1}\sin\left((2R)^{-1}d\right)}\notag\\
& \qquad\qquad\mp\frac{6m-5}{4(m-1)^24R^2d^{m-2}\sin^2\left((2R)^{-1}d\right)}\pm\frac{m^2-2m+2}{2(m-1)d^m}-\frac{(m-2)^2}{4d^2}.\tag{\theequation d}\label{sn::ex::dm::eq2d}
\end{align}
From \eqref{sn::ex::dm::eq2a}--\eqref{sn::ex::dm::eq2d}, the following series of observations can be made for various values of $m$ and sign on the cost-function depicted in \eqref{sn::ex::dm::eq1}.

\subsubsection{$(+),\ m<0$}\label{sn::ex::dm::pl}%
In this scenario, the $R$ scaling has no effect on the (A3) results. From \eqref{sn::ex::dm::eq2a}--\eqref{sn::ex::dm::eq2d}, the following is calculated:
\begin{align}
P_i(d) &\ge 0,\quad\text{for}\ i\in\{1,2\},\notag\\
P_j(d) &\le 0,\quad\text{for}\ j\in\{3,4\},\notag\\
|P_1(d)| &\ge |P_3(d)|,\label{sn::ex::dm::pl::eq1}
\end{align}
for $d\in\left(0,R\pi\right)$. Analysing limits as $y\to x$, it is observed that
\begin{align}
\lim_{d\to 0}P_k(d) &= 0,\quad\text{for}\ k\in\{1,2,3\},\notag\\
\lim_{d\to 0}P_4(d) &= -\infty.\label{sn::ex::dm::pl::eq2}
\end{align}
From \eqref{sn::ex::dm::pl::eq1} and \eqref{sn::ex::dm::pl::eq2}, it is thus concluded
\begin{equation*}
O_i(d)\ge0,\quad\text{for}\ i\in\{1,2,3\},\label{sn::ex::dm::pl::eq3}
\end{equation*}
for $d\in\left(0,R\pi\right)$; that is, the (A3) condition does not hold in this scenario.

\subsubsection{$(-),\ m<0$}\label{sn::ex::dm::nl}
Although this case bears many similarities to the scenario depicted in \subsubsecref{sn::ex::dm::pl}, it also has some notable differences. In particular, the radial scaling of the sphere does affect whether or not the (A3) condition is satisfied. However, results that are true for all $R>0$ will first be conveyed. Analysing the $P_i(d)$ expressions, it is calculated that
\begin{align}
P_i(d) &\le 0,\quad\text{for}\ i\in\{1,2\},\notag\\
P_3(d) &\ge 0,\notag\\
|P_1(d)| &\ge |P_3(d)|,\label{sn::ex::dm::nl::eq1}
\end{align}
for $d\in\left(0,R\pi\right)$ and $\forall R>0$. To proceed, the limits as $y\to x$ are calculated:
\begin{align}
\lim_{d\to 0}P_k(d) &= 0,\quad\text{for}\ k\in\{1,2,3\},\notag\\
\lim_{d\to 0}P_4(d) &= -\infty.\label{sn::ex::dm::nl::eq2}
\end{align}
From \eqref{sn::ex::dm::nl::eq1} and \eqref{sn::ex::dm::nl::eq2}, one has the following:
\begin{equation*}
O_i(d)\le0,\quad\text{for}\ i\in\{1,2,3\},
\end{equation*}
for $d\in\left(0,R\pi\right)$. So far, no mention has been made of the $O_4(d)$ expression in the current analysis. This is due to the fact that this term can switch signs on the interval $d\in\left(0,R\pi\right)$, which is a result of $P_4(d)$ becoming positive for some $d>0$. \eqref{sn::ex::dm::nl::eq2} already shows that $O_4(d)\to-\infty$ as $y\to x$; and thus, it is known that there exists some constant $M$ such that $d<M$ implies $O_4(d)<0$ by the continuity of $O_4(d)$. Thus, the analysis of the $O_4(d)$ expression is reduced to analysing the equation
\begin{align}
0 &= -2d^5\cos\left(\frac{d}{2R}\right)+(6m-5)2Rd^4\sin\left(\frac{d}{2R}\right)-(m-1)^24R^2d^3\sin^2\left(\frac{d}{2R}\right)\cos\left(\frac{d}{2R}\right)\notag\\
& \qquad\qquad-2(m-1)(m^2-2m+2)8R^3d^2\sin^3\left(\frac{d}{2R}\right)\notag\\
& \qquad\qquad-(m-1)^2(m-2)^28R^3d^m\sin^3\left(\frac{d}{2R}\right);\label{sn::ex::dm::nl::eq4}
\end{align}
the right-hand side is simply the numerator of \eqref{sn::ex::dm::eq2d}. For $m<0$ and $d\in\left(0,R\pi\right)$, the fourth term on the right-hand side of \eqref{sn::ex::dm::nl::eq4} is the only positive term. If we consider $d\in(0,1]$ and $m<0$, it is observed that
\begin{multline}
\left|2(m-1)(m^2-2m+2)R^3d^2\sin^3\left(\frac{d}{2R}\right)\right|\le
\left|(m-1)^2(m-2)^2R^3d^m\sin^3\left(\frac{d}{2R}\right)\right|,\label{sn::ex::dm::nl::eq5}
\end{multline}
as $(m-1)^2(m-2)^2>2(m-1)(m^2-2m+2)$ for all values of $m$. Thus, by \eqref{sn::ex::dm::nl::eq2} and \eqref{sn::ex::dm::nl::eq4} it is seen that $O_4(d)\le 0$ for $d\in[0,1]$. The actual value $d^*$ where $O_4(d^*)=0$ can not be explicitly represented as \eqref{sn::ex::dm::nl::eq4} is a transcendental equation. This value of $d^*$ depends on both $R$ and $m$. Even though an analytic representation of $d^*$ does not exist, it can be calculated from \eqref{sn::ex::dm::nl::eq4} that
\begin{equation}
\lim_{m\to-\infty}d^*=1,\quad\forall R>0.\label{sn::ex::dm::nl::eq6}
\end{equation} 
This can be seen from \eqref{sn::ex::dm::nl::eq4} by noticing the fourth term on the right-hand side dominates all terms for large negative values of $m$ and $d>1$ and some $R>0$ fixed. However, \eqref{sn::ex::dm::nl::eq5} holds for all $m<0$, $R>0$ and $d\le 1$; thus, one can intuitively reconcile \eqref{sn::ex::dm::nl::eq6} from \eqref{sn::ex::dm::nl::eq5} without resorting to a limit calculation. \eqref{sn::ex::dm::nl::eq6} represents the manifestation of the $R$ scaling dependency for this specific scenario. In particular, if $R\le\frac{1}{\pi}$, the (A3w) condition will be satisfied. Of course, this is not a sharp estimate; but such restrictions on the radius of the sphere will ensure the (A3w) condition is satisfied for any $m<-2$. Instead of restricting the radius of the sphere, an analogous restriction may be employed on the source and target domains: $\Diam(\sdom\cup\tdom)\le 1$. This will also ensure that the (A3w) criterion is satisfied. In order for the strong (A3) condition to be satisfied, one of the aforementioned restrictions is required plus the criterion that $\Dist(\sdom,\tdom)>0$, as
\begin{equation*}
\lim_{d\to 0}O_i(d)=0,\quad\text{for}\ i\in\{1,2,3\},
\end{equation*}
which is a straight-forward consequence of \eqref{sn::ex::dm::nl::eq2}.\\

As was already mentioned, the limit \eqref{sn::ex::dm::nl::eq6} only corresponds to a sufficient condition on $R$ for the (A3w) condition to be satisfied. The corresponding necessary condition on $R$ can indeed be calculated numerically for given values of $m$; but such a condition can not be explicitly stated due to the transcendental nature of \eqref{sn::ex::dm::nl::eq4}. However for a fixed $R$, this scenario undergoes a bifurcation in the $m$ parameter whereupon the (A3) condition will be at least weakly satisfied for all $R>0$. The value for $m$, that corresponds to this bifurcation, depends on $R$ and is again only able to be calculated numerically. One may analyse this bifurcation point as $R\to\infty$. Assuming $R>0$ and $d>1$ is such that $\frac{d}{R}<<1$, the right-hand side of \eqref{sn::ex::dm::nl::eq4} is approximated by
\begin{equation*}
(-2m^3+5m^2-4)d^5+o\left(d^7\right).
\end{equation*}
By this approximation, it is seen that if $m^*<m<0$, where $m^*$ is defined as the the negative root of the polynomial $-2m^3+5m^2-4$, then for an arbitrary fixed $d\in\left(0,R\pi\right)$
\begin{align*}
0 &> \lim_{R\to\infty}\left[-2d^5\cos\left(\frac{d}{2R}\right)+(6m-5)2Rd^4\sin\left(\frac{d}{2R}\right)\right.\\
& \qquad\qquad\left.-(m-1)^24R^2d^3\sin^2\left(\frac{d}{2R}\right)\cos\left(\frac{d}{R}\right)\right.\\
& \qquad\qquad\left.-2(m-1)(m^2-2m+2)8R^3d^2\sin^3\left(\frac{d}{2R}\right)\right].
\end{align*}
It can be calculated that
\begin{equation}
m^*\approx -.7807764064,\label{sn::ex::dm::nl::eq10}
\end{equation}
which is a numerical limit for the (A3), $m$-parameter bifurcation as $R$ tends to infinity. This is shown using elementary methods of calculus to conclude that for $m>m^*$, the first three terms of \eqref{sn::ex::dm::nl::eq4} dominate the single positive fourth term on the right-hand side of \eqref{sn::ex::dm::nl::eq4}. 
\begin{remark*}
It indeed does require the sum of all three of these terms to dominate this one positive term in the $O_4(d)$ expression; thus, there is no simple asymptotic statement to justify this behaviour outside of analysing the sign of the right-hand side of \eqref{sn::ex::dm::nl::eq4} and it's derivatives. 
\end{remark*}
Again, \eqref{sn::ex::dm::nl::eq10} correlates to a limit of $R$ tending to infinity. As $R$ decreases, this bifurcation will happen for a value of $m<m^*$. If $R$ falls below the previously discussed bifurcation point dependent on $m$ (we recall this has to be greater than $\frac{1}{\pi}$), then the (A3) condition is at least weakly satisfied. This combined with the previous statements regarding the bifurcation point in the $R$ parameter fills out the current spectrum of results regarding this scenario. The results for this particular case are summarised as follows:
\begin{itemize}
\item If $R>\frac{1}{\pi}$ and $m>m^*$, the (A3w) condition is satisfied for any $\sdom$ and $\tdom$ such that $\Cut(\sdom)\cap\tdom=\emptyset$. 
\item If $R>\frac{1}{\pi}$ and $m<m^*$, one must numerically check to see if \eqref{sn::ex::dm::nl::eq4} has real zeros for $d\in\left(0,R\pi\right)$. If this is the case, the smallest positive root of \eqref{sn::ex::dm::nl::eq4} represents an upper bound for the diameter of $\sdom\cup\tdom$ for the (A3w) condition to hold.
\item If $R\le\frac{1}{\pi}$ the (A3w) condition will hold for any $\sdom$ and $\tdom$ such that $\Cut(\sdom)\cap\tdom=\emptyset$.
\end{itemize}
All these conclusions can be strengthened to having the (A3) condition satisfied, provided $\Dist(\sdom,\tdom)>0$.

\subsubsection{$(+),\ m=0$}\label{sn::ex::dm::pe}
This, scenario is exactly the same as the case studied in \subsubsecref{sn::ex::dm::pl}, except that
\begin{equation*}
\lim_{d\to 0}P_1(d)=2,\quad\forall R>0,
\end{equation*}
which is calculated from \eqref{sn::ex::dm::eq2a}. Thus, the (A3) condition will not be satisfied in this case.

\subsubsection{$(-),\ m=0$}\label{sn::ex::dm::ne}
Straight-forward calculations from \eqref{sn::ex::dm::eq2a}--\eqref{sn::ex::dm::eq2d} indicate that
\begin{align*}
O_i(d) &< 0,\quad\text{for}\ i\in\{1,4\},\\
O_j(d) &\le 0,\quad\text{for}\ j\in\{2,3\},
\end{align*}
for $d\in\left(0,R\pi\right)$. Indeed, analysing limits, it can be ascertained that that
\begin{align*}
\lim_{d\to0}O_i(d) &= -2,\quad\text{for}\ i\in\{1,2,3\},\\
\lim_{d\to0}O_4(d) &= -\infty\\
\intertext{and}
\lim_{d\to R\pi}O_i(d) &= 0,\quad\text{for}\ i\in\{2,3\}\\
\lim_{d\to R\pi}O_j(d) &< 0,\quad\text{for}\ i\in\{1,4\}.
\end{align*}
These results are true for all values of $R>0$. Thus the (A3w) condition is satisfied independent of the radial scaling of the round sphere. To have the (A3) condition be satisfied, it is required that $\Dist(\Cut{\sdom},\tdom)>0$.

\subsubsection{$(+),\ 0<m<1$}\label{sn::ex::dm::pl1}
As with the cases where $m<0$, the scenarios for when $m>0$ also have a complex structure in that both the varying of $R$ and $m$ have bifurcations in all of the terms $O_i(d)$ for $i\in\{1,2,3,4\}$. The following exposition will be less explicit as compared to the $m<0$ cases, as the analysis is extremely similar to examples already presented in this section.\\

In this case the (A3) condition will be violated for all values of $R>0$, except in a very special situation. First, one can easily verify that $O_1(d)>0$ for $d\in\left[0,R\pi\right)$; but there are values of $0<m<1$, $R>0$ and $d$ such that $O_j(d)<0$. Specifically, for $m>m^*$ and any value of $R>0$, one has the (A3) condition being satisfied on $\mathbb{S}^2$ provided $\Dist(\sdom,\tdom)>h^*$, where $h^*$ represents the second positive root of the equation $O_4(d)=0$. $m^*$ is the root of the equation
\begin{equation}
\lim_{R\to\infty}O_4\left(R\pi\right)=0,\label{sn::ex::dm::pl1::eq1}
\end{equation}
and is approximately $.806$. If $R\lessapprox.071$, then $h^*$ is determined by the smallest positive root of $O_2(d)=0$ which is equivalent to the equation $O_3(d)=0$.
\begin{remlist*}
\item It should be noted that $O_4(d)$ has a parameter dependence on $m$ for cost-functions of the form $\scost(x,y)=\pm\frac{1}{m}d^m(x,y)$. With this, it is clear what is meant by $m^*$ being the root of \eqref{sn::ex::dm::pl1::eq1}.
\item This is one of the special cases where the dimensionality of $\mathbb{S}^n$ has an affect on the (A3) condition. Indeed, for $n>2$, one must have $O_1(d)$ less than zero, which is not the case in the current scenario. However, for $n=2$ the $O_1(d)$ is not considered, as it is clearly not possible to have a transport vector orthogonal to both arbitrary vectors $\eta$ and $\xi$ with $\eta\perp\xi$ in a two dimensional space.
\end{remlist*}

\subsubsection{$(+),\ 1<m<2$}\label{sn::ex::dm::pl2}
In this case,
\begin{equation*}
O_i(d)<0,\quad\text{for}\ i\in\{1,4\};
\end{equation*}
and $O_2(d)$ and $O_3(d)$ are monotone decreasing with
\begin{equation*}
\lim_{d\to 0}O_j(d)=\infty,\quad\text{for}\ j\in\{2,3\}.
\end{equation*}
It can be calculated that for this scenario there exists $h^*\in\left(0,R\pi\right)$ such that $O_2(h^*)=O_3(h^*)=0$. Thus, $h^*$ represents the minimal distance separation between $\sdom$ and $\tdom$ for the strong (A3) condition to be satisfied in this particular case. This conclusion is independent of radial scaling.

\subsubsection{$(+),\ 2<m<\infty$}\label{sn::ex::dm::pli}
In this case, one has that
\begin{equation*}
O_i(d)<0,\quad\text{for}\ i\in\{2,3\},
\end{equation*}
and $O_2(d)$ and $O_3(d)$ are monotone decreasing with
\begin{equation*}
\lim_{d\to 0}O_j(d)=\infty,\quad\text{for}\ j\in\{1,4\}.
\end{equation*}
The strong (A3) condition will be satisfied in this case if $\Dist(\sdom,\tdom)>h^*$, where $h^*$ is defined by the equation $O_4(h^*)=0$ (there is only one root of this equation in the interval $\left(0,R\pi\right)$). It is a straight-forward calculation to see that $h^*\in\left(0,R\pi\right)$. The results of this scenario are invariant under radial scaling.

\subsubsection{$(-),\ 0<m<1$}\label{sn::ex::dm::nl1}
Unlike the case for $m<0$, a change in sign is not necessarily tantamount to violation of the (A3) condition. Here one has that
\begin{equation*}
O_i(d)<0,\quad\text{for}\ i\in\{1,4\}
\end{equation*}
and
\begin{equation*}
O_j(0)<0,\ O_j'(d)>0,\quad\text{for}\ j\in\{2,3\}.
\end{equation*}
It is calculated that there exists a $h^*\in\left(0,R\pi\right)$ such that $O_2(h^*)=O_3(h^*)=0$. Thus, if $\Diam(\sdom\cup\tdom)<h^*$, then the (A3) condition holds.

\subsubsection{$(-),\ 1<m<2$}\label{sn::ex::dm::nl2}
In this case, it is calculated that
\begin{align*}
& O_1(d)>0,\notag\\
& O_j(0)<0,\ O_j'(d)>0,\quad\text{for}\ j\in\{2,3,4\},
\end{align*}
for all $d\in\left(0,R\pi\right)$ and $R>0$. Thus, the (A3) condition will not hold for this cost-function on $\mathbb{S}^n$ for $n>2$. On $\mathbb{S}^2$ however, there exists $h^*\in\left(0,R\pi\right)$ such that $O_4(h^*)=0$. Moreover, if $\Diam(\sdom\cup\tdom)<h^*$, then the (A3) condition will be strongly satisfied on $\mathbb{S}^2$.

\subsubsection{$(-),\ 2\le m<\infty$}\label{sn::ex::dm::nli}
Here the (A3) condition can not be satisfied for any $R>0$ with $2\le m<\infty$, as it is calculated that
\begin{equation*}
O_i(d)>0,\quad\text{for}\ i\in\{2,3\},
\end{equation*}
for all $d\in\left(0,R\pi\right)$.

\begin{remarks*}\label{sn::ex::dm::re}
This concludes a particularly long example that demonstrates the use of \lemref{sn::ac::sp::sr::lem1}. The complex structure that correlates to the above set of scenarios is due to the fact that one must analyse four different expressions for various orientations of $\xi\perp\eta$ relative to the transport vector $p$, while varying two different parameters: the exponent $m$ and the radius of the round sphere $R$. It has been demonstrated in this example that bifurcations in the (A3) behaviour happen relative to both these parameters for different orientation vectors, thus resulting in a large group of scenarios and corresponding results.
\end{remarks*}

\section{Conclusions}\label{sn::co}
The full generality of \eqref{sn::ac::eq1} has been able to be studied at the expense of utilising an inherently non-intrinsic approach in the current set of calculations. However, the use of the stereographic projection as the centrepiece of the analysis presented in this paper is uniquely powerful in the context of round spheres for two main reasons. First, the stereographic formulation is rotationally invariant on $\mathbb{S}^n$; indeed, the analysis in this paper was carried out on a arbitrary fixed $x\in\mathbb{S}^n$ against a variable Optimal Transportation target $y$. The second (and most important) simplification the stereographic formulation affords is the ability to explicitly represent the geodesic distance between two points on a sphere in terms of the projected coordinates; that is, \eqref{sn::ac::sp::hs::eq1} is valid as a representation of the geodesic distance between points $x,y\in\mathbb{S}^n$ in terms of the projected coordinates $\hat{x},\hat{y}\in\mathbb{R}^n$. Such an explicit and technically manageable representation of geodesic distance on a general Riemannian manifold is rare in any coordinate system on may choose for a chart on that manifold. Indeed, to the author's knowledge, representations of geodesic distances on even an ellipsoid result in the necessity to use highly esoteric special functions based on implicit or integral representations.\\

As the (A3) analysis in this paper reduces to the analysis of covariant derivatives of $f(d(\hat{x},\hat{y}))$, nothing can be calculated without an explicit representation of geodesic distance. Thus, given the above comments, it is computationally difficult to extend the methods in this paper to other Riemannian manifolds beyond that of the round sphere.\\

As stated in \subsecref{sn::ex::d2}, verification of the (A3) criterion is only one part to proving the regularity of potential functions associated with certain costs. The other part of proving regularity lies in the existence of gradient estimates analogous to those presented in \thmref{sn::ex::d2::thm1}, for cost-functions other than $\frac{1}{2}d^2(x,y)$. To circumvent this, the assumption that $\Cut(\sdom)\cap\tdom=\emptyset$ has been made throughout the paper. Without this assumption or a gradient estimate, the stereographic formulation becomes invalid as a point may be mapped to it's cut-locus on the sphere; and thus, move outside of the chart where the analysis was performed. In the case where $f'(d)<0$, one can not escape the assumption that $\Cut(\sdom)\cap\tdom=\emptyset$ in the stereographic formulation. In this situation, one only need to consider a case where $\Cut(\sdom)=\tdom$ to observe that the optimal mapping correlates to a mapping that takes every point in $x\in\sdom$ and maps it to that point's particular cut-locus. Thus, to analyse cases where $f'(d)<0$ free from the assumption that $\Cut(\sdom)\cap\tdom=\emptyset$, requires the use of geometrically intrinsic methods.

\bibliographystyle{amsalpha}
\bibliography{references}
\end{document}